\documentclass{article}
\usepackage{amssymb}
\usepackage{amsthm}
\usepackage{amsfonts}
\usepackage[centertags]{amsmath}

\newtheorem{theorem}{Theorem}

\newtheorem{remark}[theorem]{Remark}

\begin{document}

\date{}
\author{Thabet ABDELJAWAD\footnote{\c{C}ankaya University, Department of
Mathematics, 06530, Ankara, Turkey} , Erdal KARAPINAR\footnote{ At\i l\i m
University, Department of Mathematics, \.Incek 06836, Ankara, Turkey}}
\title{A gap in the paper "A note on cone metric fixed point theory and its equivalence"
 [Nonlinear Anal. 72(5), (2010), 2259-2261] }
\maketitle

\begin{abstract}
There is a gap in Theorem 2.2 of the paper of Du (\cite{D_2010}). In
this paper, we shall state the gap and repair it.
\end{abstract}



In 2010, Du investigated the equivalence of vectorial
versions of fixed point theorems in generalized cone metric spaces and scalar versions
of fixed point theorems in (general) metric spaces (in usual sense). He showed that
the Banach contraction principles in general metric spaces and in TVS-cone metric
spaces are equivalent. His results also extended some results of \cite{[2]} and \cite{[4]}. In this paper,
all notations are considered as in \cite{D_2010}. Further, $(E; S)$ will stand for the Hausdorff
locally convex topological vector space with S the system of seminorms generating its
topology. Also we insist on that continuity of the algebraic operations in a topological
vector space and the properties of the cone imply the relations:
\[ intP + intP \subset intP
\ \ \mbox{and}  \ \ \  \lambda intP \subset intP \ \ \mbox{for each}  \ \ \ \lambda > 0.\]
We appeal to these relations in the following. Du
proved the following result [\cite{D_2010}; Theorem 2.2].

\begin{theorem}
Let $(X,p)$ be a TVS-CMS, $x\in X$ and $\{x_n\}_{n=1}^{\infty}$ a sequence
in $X$.
Set $d_p=\xi_e\circ p$. Then the following statements hold:
\begin{itemize}
\item[($i$)] If $\{x_n\}_{n=1}^{\infty}$ converges to $x$ in TVS-CMS
$(X,p)$, then
$d_p(x_n,x)\rightarrow 0$ as $n\rightarrow \infty$,
\item[($ii$)] If $\{x_n\}_{n=1}^{\infty}$ is a Cauchy sequence in  TVS-CMS
$(X,p)$, then
$\{x_n\}_{n=1}^{\infty}$ is a Cauchy sequence (in usual sense) in $(X,d_p)$,
\item[($iii$)] If $(X,p)$ is a complete TVS-CMS, then $(X,d_p)$
is a complete metric space.
\end{itemize}
\label{lemma_eq_statements}
\end{theorem}

The author has been claimed that the conclusion $(iii)$ is immediate from conditions $(i)$ and $(ii)$.
This assertion is not true. Take a Cauchy sequence $\{x_n\}$ in $(X, d_p)$. To proceed by $(ii)$, one needs to show $\{x_n\}$ is
Cauchy in $(X, p)$, which means that the converse of the statement $(ii)$ must also hold.

In fact, the converse of  the implications of  $(i)$ and $(ii)$  hold.
We prove it here.
Regarding $(i)$ we prove that if $x_n\rightarrow x$ in $(X,d_p)$ then
$x_n\rightarrow x$ in
$(X,p)$. Let $c>>0$ be given. Take
$q \in S$ and $\delta >0$ such that $q(b)<\delta$ implies  $b<< c$.
Since $\frac{e}{n}\rightarrow 0$ in $(E,S)$, we can find $\epsilon =
\frac{1}{n_0}$ such
that $\epsilon q(e)=q(\epsilon e)<\delta $ and hence $\epsilon e << c$.
Now, choose $n_0$ such that $d_p (x_n,x)=\xi_e \circ p (x_n,x)< \epsilon$
for all $n\geq n_0$. Hence, by Lemma 1.1 $(iv)$ in \cite{D_2010}, $p
(x_n,x)<< \epsilon e<< c$ for all $n\geq n_0$. The proof of the converse
of implication $(ii)$ is similar. Now it is possible to say that $(iii)$ of
Lemma \ref{lemma_eq_statements} is immediate from the modified $(i)$ and $(ii)$.
Thus, [\cite{D_2010},Theorem 2.2] should be complete as following.

\begin{theorem}\label{summary}
Let $(X,p)$ be a TVS-CMS, $x\in X$ and $\{x_n\}_{n=1}^{\infty}$ a sequence
in $X$.
Set $d_p=\xi_e \circ p$. Then the following statements hold:
\begin{itemize}
\item[($i$)]  $\{x_n\}_{n=1}^{\infty}$ converges to $x$ in TVS-CMS
$(X,p)$ if and only if
$d_p(x_n,x)\rightarrow 0$ as $n\rightarrow \infty,$
\item[($ii$)]  $\{x_n\}_{n=1}^{\infty}$ is a Cauchy sequence in  TVS-CMS
$(X,p)$ if and only if
$\{x_n\}_{n=1}^{\infty}$ is a  Cauchy sequence  in $(X,d_p)$,
\item[($iii$)]   $(X,p)$ is a complete TVS-CMS if and only if $(X,d_p)$
is a complete metric space.
\end{itemize}
\end{theorem}
\begin{remark} \label{extra}

Above result says that for every complete TVS-cone metric
space there exists a correspondent complete usual metric space such that
the spaces are topologically isomorphic.
Actually, by using the fact that for each $c>>0$ there exist
$q \in S$ and $\delta >0$ such that $q(b)<\delta$ implies that $b<< c$, we
can show that every TVS-cone metric space is a first countable
topological  space. This is of course possible if we assume that the cone
has nonempty interior. However, there are still
many interesting fixed point theorems that could be generalized to TVS-cone
metric spaces such that their proofs do not follow directly by applying
the nonlinear scalarization function $\xi_e$ (see for example \cite{[3]}, \cite{[5]},\cite{[6]} and \cite{EK}).
\end{remark}

\textbf{Acknowledgment}
The authors express their gratitude to the referee for constructive and useful
remarks and suggestions.


\end{document}